\documentclass{amsart}
\usepackage{amsmath,amsthm,amsopn,amstext,amscd,amsfonts,amssymb,mathrsfs,mathtools}
\usepackage{dsfont}
\usepackage{comment}
\usepackage{xcolor}
\usepackage{float}
\usepackage[active]{srcltx}
\usepackage{graphicx, mathdots}

  \DeclareMathOperator{\diag}{diag}
  \DeclareMathOperator{\sgn}{sgn}
  \usepackage{mathtools}

\DeclarePairedDelimiter\floor{\lfloor}{\rfloor}

\newtheorem{theorem}{\sc Theorem}[section]
\newtheorem{conj}{\sc Conjecture}[section]

\newtheorem{proposition}{\sc Proposition}[section]

\usepackage[active]{srcltx}

\newcommand{\dps}{\displaystyle}
\usepackage{graphicx, epsfig, subfig}
\usepackage{lscape}
\usepackage{rotating}
\usepackage{caption}
\usepackage{multicol}
\usepackage{colortbl}
\usepackage{nicematrix}

\begin{document}
\title[The extreme zeros of Laguerre polynomials]{On the product of the extreme zeros of Laguerre polynomials}
\author{K. Castillo}
\address{CMUC, Department of Mathematics, University of Coimbra, 3001-501 Coimbra, Portugal}
\email{ kenier@mat.uc.pt}

\date{\today}
\begin{abstract}
The purpose of this note is twofold: firstly, it intends to bring to light  an apparently unknown property of the product of the extreme zeros of Laguerre polynomials, which in a very particular case leads to a twenty-year-old conjecture for Hermite polynomials posed by Gazeau, Josse-Michaux, and Moncea while developing numerical methods in quantum mechanics; and secondly to progress towards the solution of this problem  as an application of a parametric eigenvalue problem.
\end{abstract}
\maketitle
\section{Introduction}\label{intro}
There is an extensive literature on the properties of the zeros of orthogonal polynomials, which is rich in conjectures, remarkable theorems and proofs from The Book. This note has its motivation in one conjecture which brings to light  an apparently unknown property of the sequences of the product of the smallest and largest positive zeros of the even and odd Hermite polynomials. Hermite polynomials, $H_n$ $(n=0,1,\dots)$, are defined by the conditions
$$
\int_{-\infty}^\infty e^{-x^2} H_m(x)H_n(x)\mathrm{d}x=\sqrt{\pi}\,2^n n!\,\delta_{nm} \quad (m=0,1\dots).
$$
Let $x_{j n}$ $(j=1,\dots, \floor{n/2})$ denote the positive zeros of $H_n$ in increasing order: $
x_{1 n}<x_{2 n}<\cdots<x_{\floor{n/2} n}$.  (Since $ e^{-x^2}$ is an even function, the zeros of Hermite polynomials are situated symmetrically around the origin. For the elementary properties of the zeros of orthogonal polynomials on the real line, see \cite[Section 3.3]{S75}.) Define 
$$
y_n=x_{1 n}x_{\floor{n/2} n}\quad (n=4,5,\dots).
$$
In \cite{GMM06}\footnote{The first version of the paper was posted in arXiv in 2004.} (see also \cite{CGG10}), Gazeau,  Josse-Michaux, and Moncea, exploring what they called Hermite quantization of the real line, observe through numerical experiments that 
\begin{align}\label{inq1a}
y_n< y_{n+2},
\end{align}
or, equivalently, that $(y_{2n})_{n=2}^\infty$ and $(y_{2n+1})_{n=2}^\infty$ are monotone increasing sequences. (The reader is invited to check the validity of \eqref{inq1a} up to $n=18$ as in former times ---with paper and pencil---, using the values displayed in \cite[Section 7]{SZC52}.)  The interlacing property of the zeros of Hermite polynomials and its delicate asymptotic behavior when $n$ goes to infinity  (recall that the zeros spread out over the entire real line), quickly dispel any doubts about the complexity of the problem we are dealing with. In fact, in the event that \eqref{inq1a} were true, from the upper bounds for zeros of Hermite polynomials \cite[(6.31.19), (6.31.23)]{S75}, one sees immediately that  $(y_{2n})_{n=2}^\infty$ and $(y_{2n+1})_{n=2}^\infty$ are bounded above (see also \eqref{l1} and \eqref{l2} below), so converge, and hence these sequences are Cauchy. As if this were not enough, note that the consecutive elements of these sequences are close even for small values of $n$---  for $n=18$ and $n=17$  in \cite[Section 7]{SZC52}, it asserts, roughly, that
 \begin{align*}
 1.30382961637360&<1.32176837751291,\\[7pt]
2.58976219107561&<2.62851205461184,
\end{align*}
respectively.  
It is worth pointing out that from the convergence of $(y_{2n})_{n=2}^\infty$ and $(y_{2n+1})_{n=2}^\infty$, \cite[(6.31.19) and (6.31.23)]{S75} yields
\begin{align}
\label{l1}\lim y_{2n}&\leq\frac{\pi}{2}\lim \left\{\frac12+\frac12 \left[1-\left(\frac{\pi}{4n+1}\right)^2\right]^{1/2}\right\}^{-1/2}=\frac\pi2,\\[7pt]
\label{l2}\lim y_{2n+1}&\leq\pi\lim \left\{\frac12+\frac12 \left[1-\left(\frac{2 \pi}{4n+3}\right)^2\right]^{1/2}\right\}^{-1/2}\,\, =\pi.
\end{align}
which was also numerically observed in \cite{GMM06}. This gives a clear idea of the practical (and theoretical) importance of the inequality \eqref{inq1a}.

One fruitful approach to the zeros of Hermite polynomials  is through the use of Laguerre polynomials. Laguerre polynomials, $L^{(\alpha)}_n$ $(n=0,1,\dots; \alpha>-1)$, are defined by the conditions
 $$
\int_{0}^\infty e^{-x} x^\alpha L^{(\alpha)}_m(x)L^{(\alpha)}_n(x)\mathrm{d}x=\Gamma(\alpha+1) \binom{n+\alpha}{n}\,\delta_{nm},
$$
and Hermite polynomials can be reduced to these polynomials \cite[(5.6.1)]{S75}:
\begin{align*}
H_{2n}(x)&=(-1)^n 2^{2n} n! \,L_n^{(-1/2)}(x^2),\\[7pt]
 H_{2n+1}(x)&=(-1)^n 2^{2n+1} n!\,x\, L_n^{(1/2)}(x^2).
\end{align*}
This, in turn, invokes a more general problem: Is the sequence of the product of the extreme zeros of Laguerre polynomials, say $(y_n(\alpha))_{n=2}^\infty$,  monotone increasing on $(-1,\infty)$? In other words,  for $n=2, 3,\dots$, we are comparing the product of the extreme zeros of $L^{(\alpha)}_n$ with the product of the extreme zeros of $L^{(\alpha)}_{n+1}$. The numerical experiments suggest that it might be true. For illustration, consider $(y_n(\alpha))_{n=2}^{100}$ for some values of $\alpha$ displayed in Figure \ref{LaguerreF}. (Recall that the zeros of Laguerre are positive increasing functions of $\alpha$, see \cite[pp.122-123]{S75}.) 

\begin{figure}[H]
\centering
\includegraphics[width=10.5cm]{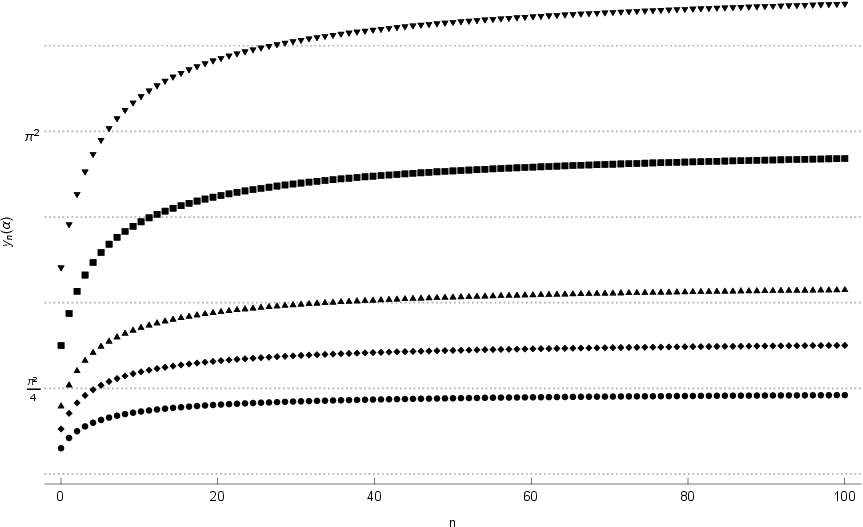}
\caption{The first $100$ elements of the sequence of the product of the extreme zeros of Laguerre polynomials for $\alpha=-1/4$, $\alpha=-1/2$, $\alpha=0$, $\alpha=1/2$, and $\alpha=1$.}
\end{figure} \label{LaguerreF}
The above question suggests an interesting (and, certainly, more complex) problem, from a physical and numerical point of view, related with the ``continuous processes" in which the extreme zeros of $L^{(\alpha)}_n$ become those of $L^{(\alpha)}_{n+1}$. In this sense we prove the next result\footnote{From \eqref{m2} below, it follows easily, for $n=5,6,\dots$, that the extreme zeros of $\widehat{L}^{(\alpha)}_{n+1}(\cdot, \alpha, 0)$ are those of $L^{(\alpha)}_{n}$ and, by definition, for all $n$, the zeros of $\widehat{L}^{(\alpha)}_{n+1}(\cdot, \alpha, 1)$ are those of $L^{(\alpha)}_{n+1}$.}.

\begin{proposition}\label{main}
Define the polynomials $\widehat{L}^{(\alpha)}_{n+1}(\cdot, t)$ $(n=0,1,\dots; \alpha>-1)$ by
$$
\widehat{L}^{(\alpha)}_{n+1}(\cdot, t)=(x-(2n+\alpha+1)){L}^{(\alpha)}_{n}(x)-n(n+\alpha)\, t^2\, {L}^{(\alpha)}_{n-1}(x),
$$
 where 
 $$
 L^{(\alpha)}_{-1}=0, \quad \widehat{L}^{(\alpha)}_{n}=(-1)^{n} n!\, L^{(\alpha)}_{n}.
 $$
Then, for all $-1<\alpha\leq 47.9603$ and $n=5,6,\dots$, the product of the extreme zeros of $\widehat{L}^{(\alpha)}_{n+1}(\cdot, t)$ is an increasing function of $t$ on $[0, 0.568774]$. \end{proposition}



What happens for the values of $\alpha>47.9603$? Unfortunately, the known estimates for the extreme zeros of the Laguerre polynomials ---the origin of which is as old as the polynomials themselves--- are too weak to allow us to conjecture anything solid for all $\alpha>-1$. After reading the following sections, the reader will be able to  obtain a better estimative depending of the values of $n$ and/or $\alpha$ fixed. For instance, it will be easy to check that the product of the extreme zeros of $\widehat{L}^{(-1/2)}_{200}(\cdot, t)$ and $\widehat{L}^{(1/2)}_{200}(\cdot, t)$ are increasing functions of $t$ on $[0, 0.8531]$.  The reader will also be able to replace the  ``transition" function $t^2$ by other suitable function (which may even depend on other parameters, for example, $\alpha$ and/or $n$).  We emphasize that the conclusion of Proposition \ref{main} does not hold on the entire interval $[0,1]$. However, after several numerical experiments, we conjecture that  the observation of Gazeau,  Josse-Michaux, and Moncea for Hermite polynomials --- $\alpha=\pm 1/2$ and $t=1$ in the present context--- is just the tip of the iceberg.

\begin{conj}\label{conj}
Assume the hypotheses and notation of Proposition \ref{main}. Then, for all $\alpha$, $n=5,6\dots$, and $t\in(0,1]$, the product of the extreme zeros of $\widehat{L}^{(\alpha)}_{n+1}(\cdot, t)$ is greater than the product of the extreme zeros of $\widehat{L}^{(\alpha)}_{n+1}(\cdot, 0)$.
\end{conj}

The proof of Proposition \ref{main} is presented in Section \ref{app} as an application of the parametric eigenvalue problem studied in Section \ref{demos}. 

\section{A parametric eigenvalue problem}\label{demos}

There is a set of real parameters, $a_{n}$ and $b_{n}>0$, associated with any (possibly finite) sequence of orthonormal polynomials on the real line, $(p_{n})_{n=0}^\infty$, 
such that
\begin{align}\label{rec}
xp_n(x)=b_{n+1}p_{n+1}(x)+a_{n+1}p_n(x)+b_np_{n-1}(x),
\end{align}
with initial conditions $p_{-1}=0$ and $p_0=1$. 
In order to recast a problem on zeros of orthogonal polynomials on the real line as a symmetric eigenvalue problem, observe from \eqref{rec} that the zeros of $p_{n+1}$ are the eigenvalues of the (unreduced) symmetric tridiagonal matrix
\begin{align}\label{m}
\mathbf{J}_{n+1}=
\begin{pmatrix}
a_1 & b_1 & & &\\[7pt]
b_1 &   a_2 &  \ddots & & \\[7pt]
& \ddots & \ddots & b_{n-1} &\\[7pt]
& & b_{n-1} & a_{n}& b_{n} \\[7pt]
&&&b_{n} & a_{n+1}
\end{pmatrix}.
\end{align}
 Let $f$ be a monotone function continuous on the close interval $[a,b]$ and differentiable on the open interval $(a, b)$ such that $f((a,b))\not=0$ and $f((a,b))\not=1$. Define the matrix-valued function
\begin{align}\label{mt}
\mathbf{J}_{n+1}(t)=\mathbf{J}_{n+1}(f, t)=
\begin{pmatrix}
a_1 & b_1 & & &\\[7pt]
b_1 &   a_2 &  \ddots & & \\[7pt]
& \ddots & \ddots &  b_{n-1} &\\[7pt]
& & b_{n-1} & a_{n}& b_n\, f(t) \\[7pt]
&&&b_n\, f(t) &  a_{n+1}
\end{pmatrix}
\end{align}
of the real variable  $t\in[a,b]$. Since $\mathbf{J}_{n+1}(t)$ is a continuous symmetric tridiagonal matrix, its eigenvalues are continuous functions of $t$ when ordered according to magnitude, $\lambda_{1}(t)\leq\lambda_{2}(t)\leq\cdots\leq\lambda_{n+1}(t)$. These inequalities are strict, with the possible exception of one of them at $t=a$ or $t=b$; that is, when $a_{n+1}$ is equal to one of the eigenvalues of $\mathbf{J}_{n}$. Note that $\mathbf{J}_{n+1}(-f, t)$ is similar to $\mathbf{J}_{n+1}(f, t)$. In fact, if we define $\mathbf{U}=\diag(1,\dots, 1, -1)$ then $\mathbf{J}_{n+1}(f, t)=\mathbf{U}\,\mathbf{J}_{n+1}(-f, t)\,\mathbf{U}$. Let $\mathbf{p}_1(t)=(p_{1 1}(t), \dots,$ $p_{n+1\, 1}(t))^\mathsf{T}$ be the unit eigenvector of $\mathbf{J}_{n+1}(t)$ associated with $\lambda_1(t)$, i.e. $\mathbf{p}_1(t)^\mathsf{T}\mathbf{p}_1(t)=1$. We thus get
 \begin{align}\label{TM}
p_{k\, 1}(t)=\frac{p_{k-1}(\lambda_1(t))}{\|\mathbf{p}_1(t)\|}\quad (k=1,2\dots, n), \quad p_{n+1\, 1}(t)=\sgn(f(t))\,\,\frac{p_{n}(\lambda_1(t))}{\|\mathbf{p}_1(t)\|},
\end{align}
where
$$
\|\mathbf{p}_1(t)\|^2=\sum_{j=0}^n p_j^2(\lambda_1(t)).
$$
From now on,  we restrict the variable $t$ to the open interval $(a, b)$. Since $\lambda_1(t)$ is a simple eigenvalue of the differentiable matrix-valued function $\mathbf{J}_{n+1}(t)$, it follows that $\lambda_1(t)$ is a differentiable function of $t$  and the Hadamard first variation formula (see \cite[(1.73)]{Tao} or, alternatively, \cite[Proposition 3.1]{CZ22}) gives
\begin{align}\label{derivada1}
\lambda'_1(t)=\mathbf{p}^\mathsf{T}_1(t)\, \mathbf{J}'_{n+1}(t)\, \mathbf{p}_1(t)=2\,b_n f'(t) p_{n 1}(t)p_{n+1\, 1}(t).
\end{align}
Let $p_{n+1}(\cdot, t)$ be defined by
\begin{align}\label{pt}
f^2(t) b_{n+1}p_{n+1}(x, t)= (x-a_{n+1}) p_n(x)-b_n f^2(t) p_{n-1}(x).
 \end{align}
We claim that the zeros of $p_{n+1}(\cdot, t)$ are the eigenvalues of $\mathbf{J}_{n+1}(t)$. Indeed, expansion in terms of the elements of the last row (or column) of the determinant $P_{n+1}(x, t)=\det(x\,\mathbf{I}_{n+1}-\mathbf{J}_{n+1}(t))$ gives
$$
P_{n+1}(x, t)= (x-a_{n+1}) P_n(x)-b^2_n f^2(t) P_{n-1}(x),
$$
$P_{n}$ and $P_{n-1}$ being,  respectively, the monic characteristic polynomials of $\mathbf{J}_{n}$ and $\mathbf{J}_{n-1}$. Since $p_n=(b_n b_{n-1}\cdots b_1)^{-1}\,P_{n}$ and $p_{n-1}=(b_{n-1}\cdots b_1)^{-1} P_{n-1}$, we have $$
p_{n+1}(\cdot, t)=(f^2(t) b_{n+1}b_n b_{n-1}\cdots b_1)^{-1}\,P_{n+1}(\cdot, t),
$$ 
which proves our assertion. Subtracting \eqref{rec} from \eqref{pt} we obtain
\begin{align}\label{relation}
f^2(t)b_{n+1}p_{n+1}(x,t)=b_{n+1}p_{n+1}(x)+b_n (1-f^2(t)) p_{n-1}(x).
\end{align}
Evaluating \eqref{relation} at $x=\lambda_1(t)$ gives
\begin{align}
\label{combination1}b_{n+1} p_{n+1}(\lambda_1(t))&=(f^2(t)-1)\,b_np_{n-1}(\lambda_1(t)).
\end{align}
 We now claim that
\begin{align}\label{norm}
 \|\mathbf{p}_1(t)\|^2=b_{n+1} p'_{n+1}(\lambda_1(t), t) p_n(\lambda_1(t))+\frac{f^2(t)-1}{f^2(t)}p^2_n(\lambda_1(t)).
\end{align}
Indeed, multiplying both sides of \eqref{pt} by $p_n(y)$ we have
\begin{align*}
f^2(t) b_{n+1}p_{n+1}(x, t)p_n(y)= (x-a_{n+1}) p_n(x)p_n(y)-b_n f^2(t) p_{n-1}(x)p_n(y).
\end{align*}
Interchanging $x$ and $y$  gives
\begin{align*}
f^2(t) b_{n+1}p_{n+1}(y, t)p_n(x)= (y-a_{n+1}) p_n(y)p_n(x)-b_n f^2(t) p_{n-1}(y)p_n(x).
\end{align*}
Now subtract both equations, and from the Christoffel–Darboux formula (see \cite[Theorem 3.2.2]{S75}) we get
\begin{align*}
&(x-y) \left(p_{n}(x)p_n(y)+f^2(t)\sum_{k=0}^{n-1} p_{k}(x)p_k(y)\right)\\[7pt]
&\quad=f^2(t) b_{n+1} (p_{n+1}(x, t)p_n(y)-p_{n+1}(y, t)p_n(x)).
\end{align*}
Hence
\begin{align*}
&p_{n}(x)p_n(y)+f^2(t)\sum_{k=0}^{n-1} p_{k}(x)p_k(y)\\[7pt]
&\quad f^2(t) b_{n+1}\left( \frac{p_{n+1}(x, t)-p_{n+1}(y,t)}{x-y} p_n(y)+ p_{n+1}(y, t)\frac{p_{n}(x)-p_{n}(y)}{x-y} \right),
\end{align*}
and  \eqref{norm} follows when $x$ tends to $y$. By the Cauchy interlacing theorem (see \cite[Exercise 1.3.14]{Tao}), $p_{n 1}(t)p_{n+1\,1}(t)\not=0$. We also claim that $\lambda_1(t)-a_{n+1}\not=0$. To obtain a contradiction,  suppose that $\lambda_1(t)=a_{n+1}$, and so $p_{n-1}(\lambda_1(t))=0$. Using \eqref{norm} and \eqref{pt} we obtain
\begin{align}
\label{paige} \frac{1}{p_{n 1}(t)p_{n+1\,1}(t)}&=\sgn(f(t))\left(\frac{b_{n+1} p'_{n+1}(\lambda_1(t), t)}{p_{n-1}(\lambda_1(t))}+(f^2(t)-1)\frac{b_n}{\lambda_1(t)-a_{n+1}}\right)\\[7pt]
\label{paige2}&=\sgn(f(t))\,\frac{p_{n}(\lambda_1(t))}{p_{n-1}(\lambda_1(t))}\frac{1}{p^2_{n+1\, 1}(t)}.
\end{align}
 (For abbreviation, we let $p'_{n+1}(\cdot, t)$ stand for the derivative with respect to the indeterminate of the polynomial.) Combining \eqref{combination1} and \eqref{paige} we can rewrite \eqref{derivada1} as
\begin{align}\label{C1}
\frac{1}{\lambda'_{1}(t)}=\sgn(f(t))\,\frac{f^2(t)-1}{2\,f'(t)} \left(\frac{p'_{n+1}(\lambda_{1}(t), t)}{p_{n+1}(\lambda_1(t))}+\frac{1}{\lambda_1(t)-a_{n+1}}\right).
\end{align}
Moreover, combining \eqref{derivada1} and \eqref{paige2}, and using the Cauchy interlacing theorem, we have 
\begin{align}\label{de1}
\sgn(\lambda'_{1}(t))=-\sgn(f(t) f'(t)).
\end{align} 
From \eqref{C1}, taking the partial derivative of \eqref{relation} with respect to $x$ and using \eqref{combination1}, we obtain
\begin{align}
\label{F} \frac{\lambda_{1}(t)}{\lambda'_{1}(t)}=\sgn(f(t))\,\frac{f^2(t)-1}{2f^2(t) f'(t)} \left( P(\lambda_1(t))+\frac{\lambda_1(t) f^2(t)}{\lambda_1(t)-a_{n+1}}\right),
\end{align}
where
$$
P(x)=x\,\frac{p'_{n+1}(x)}{p_{n+1}(x)}-x\,\frac{p'_{n-1}(x)}{p_{n-1}(x)}.
$$
By the same arguments, we also have 
\begin{align*}
\frac{\lambda_{n+1}(t)}{\lambda'_{n+1}(t)}=\sgn(f(t))\,\frac{f^2(t)-1}{2 f^2(t) f'(t)} \left( P(\lambda_{n+1}(t))+\frac{\lambda_{n+1}(t) f^2(t)}{\lambda_{n+1}(t)-a_{n+1}}\right)
\end{align*} 
and 
\begin{align}\label{de2}
\sgn(\lambda'_{n+1}(t))=\sgn(f(t) f'(t)).
\end{align}
Now note that $\lambda_1(t)\lambda_{n+1}(t)$ is continuous on the close interval $[a,b]$ and differentiable on the open interval $(a, b)$. 
By the Mean Value Theorem, there is a point $c$ (strictly) between $a$ and $b$ for which
\begin{align*}
&\lambda_1(b)\lambda_{n+1}(b)-\lambda_1(a)\lambda_{n+1}(a)\\[7pt]
&\quad =(b-a)\lambda'_1(c)\lambda'_{n+1}(c)\left(  \frac{\lambda_{1}(c)}{\lambda'_{1}(c)}+\frac{\lambda_{n+1}(c)}{\lambda'_{n+1}(c)}\right)\\[7pt]
&\quad =\sgn(f(c))\,\frac{f^2(c)-1}{2 f^2(c) f'(c)} \,(b-a)\lambda'_1(c)\lambda'_{n+1}(c)\\[7pt]
&\qquad \times \left(P(\lambda_{1}(c))+P(\lambda_{n+1}(c))+f^2(c)\left(\frac{\lambda_1(c) }{\lambda_1(c)-a_{n+1}(\alpha)}+\frac{\lambda_{n+1}(c)}{\lambda_{n+1}(c)-a_{n+1}(\alpha)} \right)\right).
\end{align*}
Finally, under our hypotheses, 
$$
\sgn\left(\sgn(f(c))\lambda'_1(c)\lambda'_{n+1}(c) \frac{f^2(c)-1}{f^2(t) f'(c)} \right)=-\sgn(f(c)f'(c))\sgn(f^2(c)-1),
$$
and we are thus led to the following result.

\begin{theorem}\label{lemmamain}
Let $f$ be a monotone function continuous on the close interval $[a,b]$ and differentiable on the open interval $(a, b)$ such that $f((a,b))\not=0$ and $f((a,b))\not=1$. Define the matrix-valued function
\begin{align*}
\mathbf{J}_{n+1}(t)=
\begin{pmatrix}
a_1 & b_1 & & &\\[7pt]
b_1 &   a_2 &  \ddots & & \\[7pt]
& \ddots & \ddots &  b_{n-1} &\\[7pt]
& & b_{n-1} & a_{n}& b_n\, f(t) \\[7pt]
&&&b_n\, f(t) &  a_{n+1}
\end{pmatrix}
\end{align*}
of the real variable  $t\in[a,b]$.
Let $\lambda_{1}(t)$ and $\lambda_{n+1}(t)$ be the extreme eigenvalues of $\mathbf{J}_{n+1}(t)$. Define
$$
P(x)=x\,\frac{P'_{n+1}(x)}{P_{n+1}(x)}-x\, \frac{P'_{n-1}(x)}{P_{n-1}(x)},
$$
where $P_{n+1}$ and $P_{n-1}$ are the characteristic polynomials of $\mathbf{J}_{n+1}(b)$ and its $(n-1)$-by-$(n-1)$ principal submatrix, respectively. Then there is a point $c\in(a, b)$ for which
\begin{align*}
&\sgn\left(\lambda_1(b)\lambda_{n+1}(b)-\lambda_1(a)\lambda_{n+1}(a)\right)=-\frac{\sgn(f(c)f'(c))}{\sgn(f^2(c)-1)}\\[7pt]
&\quad \times \sgn\left(P(\lambda_{1}(c))+P(\lambda_{n+1}(c))+f^2(c)\left(\frac{\lambda_1(c) }{\lambda_1(c)-a_{n+1}}+\frac{\lambda_{n+1}(c)}{\lambda_{n+1}(c)-a_{n+1}} \right)\right).
\end{align*}
\end{theorem}

\section{Proof of Proposition \ref{main} and further results}\label{app}
From the recurrence formula for Laguerre polynomials \cite[(5.1.10)]{S75}, it is easy to see that the zeros of the (orthonormal) polynomials
\begin{align*}
p_{n+1}=(-1)^{(n+1)} \binom{n+\alpha+1}{n+1}^{-1/2}\, L^{(\alpha)}_{n+1}
\end{align*}
$(n=0,1,\dots; \alpha>-1)$ are the eigenvalues of the symmetric tridiagonal matrix $\mathbf{J}_{n+1}=\mathbf{J}_{n+1}(\alpha)$ given by \eqref{m} with\begin{align*}
a_{n+1}&=a_{n+1}(\alpha)=2n+\alpha+1,\\[7pt]
b_n&=b_n(\alpha)=\sqrt{n(n+\alpha)}.
\end{align*}
Using these values, define the matrix-valued function $\mathbf{J}_{n+1}(t)=\mathbf{J}_{n+1}(t, \alpha)$ as in Theorem \ref{lemmamain}, i.e. 
$$
\mathbf{J}_{n+1}(t, \alpha)=\begin{pmatrix}
\alpha+1 & \sqrt{\alpha+1} & &\\[7pt]
\ \ \sqrt{\alpha+1} &   \alpha+3 &  \ddots & \\[7pt]
& \ddots & \ddots &  \sqrt{n(n+\alpha)}\, f(t)\\[7pt]
& & \sqrt{n(n+\alpha)}\, f(t)& 2n+\alpha+1
\end{pmatrix},$$
for an appropriate function $f$ to be defined later. Let $\lambda_{1}(t)$ and $\lambda_{n+1}(t)$ be the extreme eigenvalues of $\mathbf{J}_{n+1}(t, \alpha)$. From \eqref{combination1} we have
\begin{align}\label{rel1}
L_{n+1}^{(\alpha)}(\lambda_1(c))=\frac{(n+\alpha)(f^2(c)-1)}{n+1}   L_{n-1}^{(\alpha)}(\lambda_1(c)),
\end{align}
for all $c\in(0,1)$. Similarly, \eqref{pt} yields
\begin{align}\label{rel2}
L_{n}^{(\alpha)}(\lambda_1(c))=-\frac{(n+\alpha)\, f^2(c)}{\lambda_1(c)-a_{n+1}}L_{n-1}^{(\alpha)}(\lambda_1(c)).
\end{align}
(Since the Laguerre polynomials have exclusively positive zeros \cite[Theorem 3.3.1]{S75}, the $n$-by-$n$ principal submatrix of $\mathbf{J}_{n+1}(t, \alpha)$, $\mathbf{J}_{n}(\alpha)$, is an oscillatory matrix \cite[Theorem 11, p. 103]{GK02}, and the majorization theorem \cite[p. 182]{Ga82} implies that $\lambda_{1}(0)<a_n$\footnote{Of course, in this particular case, we can also use a bound proved in 1933 by Hahn. However, this and other bounds are consequences of more general properties. For instance, the majorization theorem for oscillatory matrices also implies that $2n+\alpha-1$ is a lower bound for the largest eigenvalue of $\mathbf{J}_{n}(\alpha)$. This rough estimate was proved by Szeg\H{o} \cite[(6.2.14)]{S75} from a theorem due to Laguerre.}. Therefore $\lambda_1(c)-a_{n+1}<0$, because $\lambda_1(0)-a_{n+1}<0$ and $\lambda_1(t)$ is a decreasing function of $t$ on $(0,1)$.) Now, recall that \cite[(5.1.13), (5.1.14)]{S75}
\begin{align}
\label{P1} x  \frac{\mathrm{d}}{\mathrm{d}x}L^{(\alpha)}_{n}(x)&=n L^{(\alpha)}_{n}(x)-(n+\alpha) L^{(\alpha)}_{n-1}(x),\\[7pt]
\label{P2} \frac{\mathrm{d}}{\mathrm{d}x}L^{(\alpha)}_{n-1}(x)&=\frac{\mathrm{d}}{\mathrm{d}x}L^{(\alpha)}_{n}(x)+L^{(\alpha)}_{n-1}(x).
\end{align}
Since $c\in(0,1)$, we have $L^{(\alpha)}_{n+1}(\lambda_{1}(c))\not=0$. From \eqref{P1}, \eqref{rel1}, and \eqref{rel2} we deduce that
\begin{align*}
\frac{\lambda_{1}(c)\, \dps\frac{\mathrm{d}}{\mathrm{d}x}L^{(\alpha)}_{n+1}(\lambda_{1}(c))}{L^{(\alpha)}_{n+1}(\lambda_{1}(c))}
=n+1+b_{n+1}^2\,\frac{f^2(c)}{f^2(c)-1}\frac{1}{\lambda_1(c)-a_{n+1}}.
\end{align*}
Likewise, from \eqref{P2}, \eqref{P1}, and \eqref{rel2}  we get
\begin{align*}
\frac{\lambda_{1}(c)\, \dps\frac{\mathrm{d}}{\mathrm{d}x}L^{(\alpha)}_{n-1}(\lambda_{1}(c))}{L^{(\alpha)}_{n-1}(\lambda_{1}(c))}
=\lambda_1(c)-n-\alpha-b_n^2\, f^2(c) \frac{1}{\lambda_1(c)-a_{n+1}}.
\end{align*}
It follows from what has already been proved that 
\begin{align*}
P(\lambda_1(c))&=a_{n+1}-\lambda_1(c)+(a_{n+1}+b_n^2\, f^2(c))\,\frac{f^2(c)}{f^2(c)-1}\,\frac{1}{\lambda_1(c)-a_{n+1}}.
\end{align*}
And the same applies, mutatis mutandis, for $P(\lambda_{n+1}(c))$.  A simple computation gives
\begin{align}\label{aux}
\nonumber&Q(c)=P(\lambda_1(c))+P(\lambda_{n+1}(c))+f^2(c)\left(\frac{\lambda_1(c) }{\lambda_1(c)-a_{n+1}}+\frac{\lambda_{n+1}(c)}{\lambda_{n+1}(c)-a_{n+1}} \right)\\[7pt]
 &\quad=(2 a_{n+1}-\lambda_1(c)-\lambda_{n+1}(c))\\[7pt]
\nonumber &\quad\times\left(1-\frac{f^4(c)}{1-f^2(c)}\, \frac{a_{n+1}+b_n^2}{(a_{n+1}-\lambda_1(c))(\lambda_{n+1}(c)-a_{n+1})}\right)+2f^2(c).
\end{align}
Set $f(t)=t$ on $[0, 1/2]$. (As we will see below, we could replace $1/2$ with $0.568774$. The details are left to the reader.) Since $\sgn(f(t)f'(t))>0$, \eqref{de1} and \eqref{de2} imply, respectively, $\lambda'_1(t)<0$ and $\lambda'_{n+1}(t)>0$. Let $\lambda_{1\, n+1}$ and $\lambda_{n+1\, n+1}$ be the extreme eigenvalues of $\mathbf{J}_{n+1}(\alpha)$.  By \cite[Theorem 1]{DK10}  we also have
\begin{align}\label{dk}
\lambda_{n+1\, n+1}&\leq\frac{2(n+1)^2+2(\alpha+1)+(\alpha-1)(n+1)+2n\sqrt{(n+1)^2+(\alpha+1)(n+3)}}{n+3},
\end{align}
for $\alpha>-1$. By \cite[(5) and (7)]{DJ12}\footnote{There are some typos in the information of both expressions there.} we have 
\begin{align}\label{m0}
\nonumber&\lambda_{1\, n}<\frac{(\alpha+2)_2(3n+2\alpha+2)}{2(n+\alpha+1)_2}\\[7pt]
 &-\frac{\dps\sqrt{(\alpha+2)_2(9(\alpha+2)_2+2(2\alpha+5)(\alpha^2+5\alpha+10)(n-1)+(5\alpha^2+25\alpha+38)(n-1)^2)}}{2(n+\alpha+1)_2},
\end{align}
and
\begin{align}\label{m2}
2n+\alpha-2+\sqrt{n^2-2n+\alpha\,n+2}<\lambda_{n\, n}.
\end{align}
 From \eqref{dk}, \eqref{m0}, and the monotonicity of $\lambda_1(t)$ and $\lambda_{n+1}(t)$ we get
\begin{align}\label{prod1}
\nonumber&2 a_{n+1}-\lambda_{1}(c)-\lambda_{n+1}(c)>2 a_{n+1}-\lambda_{1}(0)-\lambda_{n+1}(1/2)>2 a_{n+1}-\lambda_{1\, n}-\lambda_{n+1\, n+1}\\[7pt]
\nonumber&>4n+2\alpha+2-\frac{(\alpha+2)_2(3n+2\alpha+2)}{2(n+\alpha+1)_2}\\[7pt]
\nonumber&+\frac{\dps\sqrt{(\alpha+2)_2(9(\alpha+2)_2+2(2\alpha+5)(\alpha^2+5\alpha+10)(n-1)+(5\alpha^2+25\alpha+38)(n-1)^2)}}{2(n+\alpha+1)_2}\\[7pt]
&-\frac{2(n+1)^2+2(\alpha+1)+(\alpha-1)(n+1)+2n\sqrt{(n+1)^2+(\alpha+1)(n+3)}}{n+3}>0,
\end{align}
for $n=5, 6,\dots$ and $-1<\alpha\leq 47.9603$. Using \eqref{m0}, \eqref{m2}, and the monotonicity of $\lambda_{1}(t)$ and $\lambda_{n+1}(t)$, we deduce that
 \begin{align}\label{prod2}
&\nonumber(a_{n+1}-\lambda_1(c))(\lambda_{n+1}(c)-a_{n+1})>(a_{n+1}-\lambda_{1\, n})(\lambda_{n\, n}-a_{n+1})\\[7pt]
& \nonumber>\left(2n+\alpha+1-\frac{(\alpha+2)_2(3n+2\alpha+2)}{2(n+\alpha+1)_2}\right.\\[7pt]
&\nonumber+\left.\frac{\dps\sqrt{(\alpha+2)_2(9(\alpha+2)_2+2(2\alpha+5)(\alpha^2+5\alpha+10)(n-1)+(5\alpha^2+25\alpha+38)(n-1)^2)}}{2(n+\alpha+1)_2}\right)\\[7pt]
&\nonumber\times \left(\sqrt{n^2-2n+\alpha\,n+2}-3\right)>\frac{1}{12} (n(\alpha+2)+n^2+\alpha+1)=\frac{1}{12} \left(a_{n+1}+b_n^2\right)\\[7pt]
&>\frac{c^4}{1-c^2}\,(a_{n+1}+b_n^2)=\frac{f^4(c)}{1-f^2(c)}\,(a_{n+1}+b_n^2),
\end{align}
for $n=5,6,\dots$ and $\alpha>-1$, and so $Q(c)>0$  for $f(t)=t$, $-1<\alpha\leq 47.9603$, and $c\in(0, 1/2)$. By Theorem \ref{lemmamain}, 
$$
\lambda_1(1/2)\lambda_{n+1}(1/2)>\lambda_1(t)\lambda_{n+1}(t)>\lambda_1(0)\lambda_{n+1}(0)
$$
 for $f(t)=t$, $n=5,6,\dots$, $-1<\alpha\leq 47.9603$, and $t\in(0, 1/2)$. If we return to the inequality \eqref{prod2}, putting
 $``0.985135"$ instead of $``1/12\approx 0.0833333"$ we can also conclude, for instance, that 
\begin{align}\label{inq1}
\lambda_1(0.784522)\lambda_{n+1}(0.784522)>\lambda_1(0)\lambda_{n+1}(0),
\end{align}
for $f(t)=t$, $n=100, 101,\dots$, and $-1<\alpha\leq 47.9603$. Indeed, what we have is that $\lambda_1(t)\lambda_{n+1}(t)$ is an increasing function of $t$ on $[0, 0.784522]$ for these values of $n$ and $\alpha$. Although even this result can be slightly improved, the reader is warned that $\lambda_1(t)\lambda_{n+1}(t)$ is not monotonic on $[0,1]$. For illustration, consider $\lambda_1(t)\lambda_{n+1}(t)$  for some function $f(t)$, $\alpha=1/2$, and $n=20$ displayed in Figure \ref{Prod}.
\begin{figure}[H]
\centering
\includegraphics[width=10.5cm]{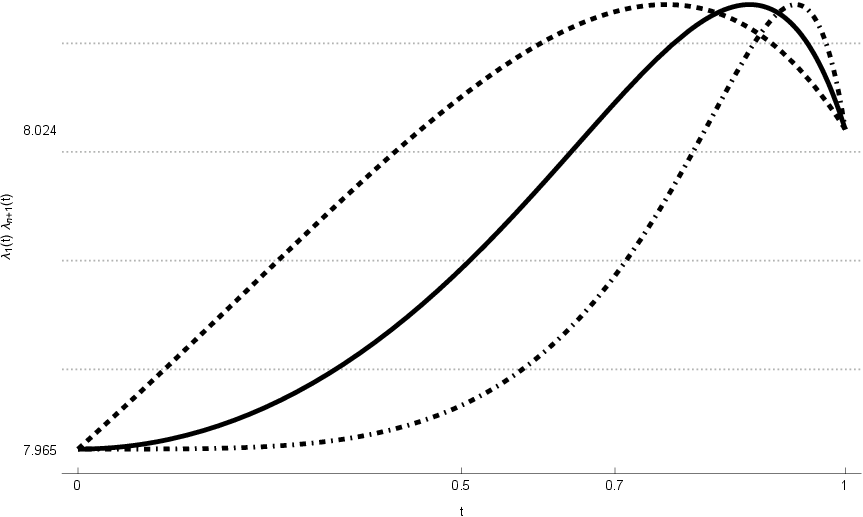}
\caption{$\lambda_1(t)\lambda_{n+1}(t)$ for $f(t)=t$ (solid line), $f(t)=\sqrt{t}$ (dashed line), $f(t)=t^2$ (dash-dotted line), $\alpha=1/2$, and $n=20$.}
\end{figure} \label{Prod}
Finally, note that for all $t$, the eigenvalues of $\mathbf{J}_{n+1}(t, \alpha)$ (and, in particular, the product of its extreme eigenvalues) are increasing functions of $\alpha$. Indeed,
$$
\frac{\mathrm{d}\mathbf{J}_{n+1}}{\mathrm{d}\alpha}(t, \alpha)=\begin{pmatrix}
1 & \dps\frac{1}{2\,\sqrt{\alpha+1}} & &\\[7pt]
\ \ \dps\frac{1}{2\,\sqrt{\alpha+1}}&   1 &  \ddots & \\[7pt]
& \ddots & \ddots & \dps\frac{n\, t}{2\,\sqrt{n(n+\alpha)}}\\[7pt]
& & \dps\frac{n\, t}{2\,\sqrt{n(n+\alpha)}} & 1
\end{pmatrix},$$
and so $(\mathrm{d}\,\mathbf{J}_{n+1})(\mathrm{d}\,\alpha)(t, \alpha)$ is strictly diagonal dominant with positive diagonal elements. From this we see that $(\mathrm{d}\,\mathbf{J}_{n+1})(\mathrm{d}\,\alpha)(t, \alpha)$ is positive definite \cite[Corollary 7.2.3]{HJ} and the result follows from the Hadamard first variation formula.

\section*{Acknowledgements}
The author is indebted to F. R. Rafaeli for drawing his attention to the work \cite{GMM06}. This work is partially supported by the Centre for Mathematics of the University of Coimbra (funded by the Portuguese Government through FCT/MCTES, DOI 10.54499/UIDB/00324/2020). The author is also supported by FCT through DOI 10.54499/2022.00143.CEECIND/CP1714/CT0002. 
\bibliographystyle{plain}
\bibliography{bib} 

\begin{thebibliography}{10}

\bibitem{CZ22}
K.~Casillo and I.~Zaballa.
\newblock On a formula of {T}hompson and {M}c{E}nteggert for the adjugate
  matrix.
\newblock {\em Linear Algebra Appl.}, 634:37--56, 2022.

\bibitem{CGG10}
N.~Cotfas, J.~P. Gazeau, and K.~G{\'o}rska.
\newblock Complex and real {H}ermite polynomials and related quantizations.
\newblock {\em J. Phys. A: Math. Theor.}, page 305304, 2010.

\bibitem{DK10}
D.~K. Dimitrov and G.~P. Nikolov.
\newblock Sharp bounds for the extreme zeros of classical orthogonal
  polynomials.
\newblock {\em J. Approx. Theory}, 162:1793--1804, 2010.

\bibitem{DJ12}
K.~Driver and K.~Jordaan.
\newblock Bounds for extreme zeros of some classical orthogonal polynomials.
\newblock {\em J. Approx. Theory}, 164:1200--1204, 2012.

\bibitem{GK02}
F.~P. Gantmacher and M.~G. Krein.
\newblock {\em Oscillation matrices and kernels and small vibrations of
  mechanical systems. (Translation based on the 1941 {R}ussian original.
  {E}dited and with a preface by {A}lex {E}remenko)}.
\newblock AMS Chelsea Publishing, Providence, RI, revised edition, 2002.

\bibitem{Ga82}
J.~Garloff.
\newblock Majorization between the diagonal elements and the eigenvalues of an
  oscillating matrix.
\newblock {\em Linear Algebra Appl.}, 47:181--184, 1982.

\bibitem{HJ}
R.~A. Horn and C.~R. Johnson.
\newblock {\em Matrix Analysis}.
\newblock Cambridge University Press, New York, second edition, 2013.

\bibitem{GMM06}
J.~P. Gazeau F.~X. Michaux and P.~Monceau.
\newblock Finite dimensional quantizations of the (q,p) plane: new space and
  momentum (or quadratures) inequalities.
\newblock {\em Modern Phys. B}, 20:1778--1791, 2006.

\bibitem{SZC52}
H.~E. Salzer, R.~Zucker, and R.~Capuano.
\newblock Table of the zeros and weight factors of the first twenty {H}ermite
  polynomials.
\newblock {\em J. Research Nat. Bur. Standards}, 48:111--116, 1952.

\bibitem{S75}
G.~Szeg\H{o}.
\newblock {\em Orthogonal polynomials}, volume~23.
\newblock Amer. Math. Soc. Coll. Publ., Amer. Math. Soc., Providence, R. I.,
  {F}ourth edition, 1975.

\bibitem{Tao}
T.~Tao.
\newblock {\em Topics in random matrix theory}, volume 132 of {\em Grad. Stud.
  Math.}
\newblock American Mathematical Society, Providence, RI, 2012.

\end{thebibliography}
 \end{document}